\documentclass{amsart}
\usepackage{graphicx}
\usepackage{latexsym}
\usepackage{amsmath}
\usepackage{amssymb}
\usepackage{amsthm}
\usepackage{amscd}
\usepackage{amsrefs}

\newtheorem{thm}{Theorem}[section]
\newtheorem{lemma}[thm]{Lemma}
\newtheorem{prop}[thm]{Proposition}
\newtheorem{cor}[thm]{Corollary}

\theoremstyle{definition}
\newtheorem{definition}[thm]{Definition}
\newtheorem{remark}[thm]{Remark}
\newtheorem{ex}[thm]{Example}

\newcommand{\al}{\alpha}                
\newcommand{\sig}{\sigma}               

\newcommand{\om}{\omega}                

\newcommand{\ra}{\rightarrow}           
\newcommand{\lra}{\longrightarrow}

\newcommand{\subs}{\subseteq}           

\newcommand{\sect}{{\mathcal{x}}}
\newcommand{\e}{{\acute{\textup{e}}}}
\newcommand{\et}{{\acute{e}t}}
\newcommand{\Z}{{\mathbb{Z}}}
\newcommand{\Ze}{{\Z_{\et}}}

\newcommand{\Q}{{\mathbb{Q}}}

\newcommand{\Hc}{{\mathcal{H}}}
\newcommand{\HT}{{\widehat{H}}}

\newcommand{\trg}{{\tau_{_{\leq 0}} \textup{R}\Gamma}}
\newcommand{\RG}{{\textup{R}\Gamma}}
\newcommand{\tn}{{\tau_{_{\leq 0}}}}

\newcommand{\ag}{{\alpha_{L/K}}}
\newcommand{\bg}{{\beta_{L/K}}}
\newcommand{\ah}{{\alpha_{L/E}}}
\newcommand{\bh}{{\beta_{L/E}}}

\newcommand{\inv}{{\textup{inv}}}
\newcommand{\Gal}{{\textup{Gal}}}
\newcommand{\Hom}{{\textup{Hom}}}

\newcommand{\coker}{{\textup{coker}}}

\newcommand{\infl}{{\textup{Infl}}}

\newcommand{\res}{{\textup{Res}}}
\newcommand{\Res}{{\textup{Res}}}
\newcommand{\image}{{\textup{Im}}}
\newcommand{\W}{{\textup{Weil}}}
\newcommand{\U}{{\textup{U}}}

\newcommand{\incl}{{\textup{incl.}}}

\author{Karen Acquista}
\title{The Weil group of a hyper-class formation}
\date{\today}
\email{kea@math.bu.edu}
\address{Boston University \\ Department of Mathematics and Statistics \\ Boston, MA}

\begin{document}

\begin{abstract}  In the setting of generalized class field theory, one has a field $K$ with a hyper-class formation $C^{\bullet}$, a complex of $G_K$-modules that plays the role of a class formation.  Letting $C(L)$ denote $H^0(G_L, C^{\bullet})$, there is a distinguished cohomology class in $H^2(\Gal(L/K), C(L))$, the hyper-fundamental class.  In this paper, we prove that a group extension of $\Gal(L/K)$ by $C(L)$ whose corresponding cohomology class is the hyper-fundamental class satisfies the classical axioms of a Weil group.
\end{abstract}

\maketitle


Classically, a Weil group is associated to a field $K$ with a class formation $C$.  
For any finite Galois extension $L/K$, 
it is constructed 
as a group extension of the Galois group $\Gal(L/K)$ by the formation module $C^{G_L}$, whose corresponding cohomology class is 
the canonical generator of $H^2(\Gal(L/K), C^{G_L})$.  That is, the Weil group 
is a ``twisted product'' of the Galois group and the formation module, with a multiplication law defined using a 2-cocycle representing the fundamental class.

The Weil group of a local or global field appears in a variety of arithmetic settings,
see Tate \cite{TateNT}; 
most recently,  Lichtenbaum \cite{Li2} has used the Weil group to construct a new cohomology theory for global fields that appears to capture integral arithmetic information.  Weil-$\e$tale cohomology groups are expected to be finitely generated over $\Z$, and contain {\em more} information than their $\e$tale counterparts.  Also, special values of the zeta functions appear to be genuine Euler products in this cohomology theory.  The main drawback of the Weil-$\e$tale cohomology theory is that it is only defined when there is a Weil group, that is, in the comparatively narrow setting of class field theory.  

In the 1980's, Bloch, Kato, Saito, and Parshin obtained many new results about abelian extensions 
of fields 
beyond the classical scope of local and global fields, 
see Raskind \cite{Ras} for a survey.  In the setting of generalized class field theory in which we will place ourselves, we are given a field $K$ with a {hyper}-class formation $C^{\bullet}$, a complex of $G_K$-modules that plays the role of a class formation.  The purpose of this note is to show that given a hyper-class formation, one can construct a Weil group for any finite Galois extension, in the classical sense.

In this exposition, we first review the theory of hyper-class formations as was introduced by Koya \cite{Ko}, and prove some results about the Tate cohomology groups $\HT^q(G, A^{\bullet})$.  Here, $G$ is a finite group, and $A^{\bullet}$ is bounded complex of $G$-modules.  Using the classical tool of dimension-shifting, we obtain a surprising general result: there exists a $G$-module $A$ such that for any subgroup $H \subs G$ and any $q \in \Z$, 
$$
\HT^q(H, A^{\bullet}) \simeq \HT^q(H, A),
$$
see Proposition \ref{cshift}.  This result can be used to obtain a new proof of the Tate-Nakayama Theorem for complexes (Theorem \ref{tn}); we also easily see, as in the classical setting, that a hyper-class formation gives rise to a description of the abelianized Galois group.

Next, we include a treatment of the axiomatic construction of the Weil group of a class formation, along the lines of Chapter 14 of Artin-Tate \cite{AT}.  
Then, we construct a Weil group associated to a field $K$ with a hyper-class formation $C^{\bullet}$.  Letting $C(L)$ denote $H^0(G_L, C^{\bullet})$, there is a distinguished element of $H^2(\Gal(L/K), C(L))$, 
the image of the fundamental class, 
see Definition \ref{hfun}.  In Theorem \ref{hweil}, we prove that a group extension of $\Gal(L/K)$ by $C(L)$ corresponding to 
this distinguished class
is a Weil group for $L/K$; dimension-shifting for complexes plays a key role in the proof.  As in the classical setting, the ``Artin map'' coming from the hyper-class formation can be recovered by abelianizing the Weil group.

Lastly, we discuss some examples of hyper-class formations, and their associated Weil groups.  The first new example is associated to a 2-local field, a complete discrete valuation field with residue field a local field.  Koya \cite{Ko} 
showed that Lichtenbaum's weight two motivic complex is a hyper-class formation.  We explain the relationship between a Weil group for a 2-local field and a Weil group for its residue field.  We 
have similar 
results for the Weil group associated to an unramified extension of a complete discrete valuation field with residue field a global field.

The best description of the Weil group of a local or global field can only be obtained after passing over a limit of all finite extensions.  
Conditions are given in Chapter 14 of Artin-Tate \cite{AT} 
ensuring the existence of a compatible system of Weil groups associated to a topological class formation.
Unfortunately, these nontrivial topological considerations prevent us from giving such a construction for our new examples of Weil groups.  If $K$ is a 2-local field, there does not appear to be a natural way to make $K_2(K)$ into a locally compact topological group, see Fesenko \cite{F} for 
one approach.
For the purposes of arithmetic applications, it will be crucial to find a way to pass to the limit. 

{\sc Notations:} 
In this exposition, we place ourselves in the derived category whenever it is appropriate.  This gives rise to minor differences in the statements of certain Definitions and Theorems, compared to Koya \cite{Ko}.  

Let $G$ be a group.  The derived category of $G$-modules is the category of complexes of $G$-modules, up to quasi-isomorphism.  If $A^{\bullet}$, then the homology of $A^{\bullet}$, denoted $\Hc^q(A^{\bullet})$, is an invariant of the class of $A^{\bullet}$ in the derived category.  Given an element $A^{\bullet}$ of the derived category of $G$-modules, we say that $A^{\bullet}$ is {\em bounded below} if $\Hc^q(A^{\bullet})$ vanishes for $q$ small enough.
Let $A^{\bullet}$ be a bounded below element of the derived category of $G$-modules.  An {\em injective resolution} is a complex $I^{\bullet}$ of injective $G$-modules that represents $A^{\bullet}$.  Given any bounded below complex $A^{\bullet}$, one can find a complex $I^{\bullet}$ of injectives and a quasi-isomorphism $A^{\bullet} \lra I^{\bullet}$; this will also be called an injective resolution.  Injective resolutions are used to compute right derived functors.  
The truncation functor, denoted $\tau_{_{\leq n}}$, is compatible with quasi-isomorphisms.  That is, truncation is well-defined in the derived category.    
The functor of taking $G$-invariants will be denoted $\Gamma(G, -)$; its derived functor will be denoted $\RG(G, -)$.  It has the property that:
$$
\Hc^q(\RG(G, A^{\bullet})) = H^q(G, A^{\bullet}),
$$
for all $q \in \Z$.  These notations are compatible with Weibel, see \cite{W} for details.  

\section{Hyper-class formations}

Let $K$ be a field, and $G_K$ its absolute Galois group.  The following definition was 
introduced by Koya \cite{Ko}:

\begin{definition} \label{hclassformation}  Let $C^{\bullet}$ be a 
bounded below element in the derived category of 
$G_K$-modules, acyclic in positive degrees.  Suppose that for every finite Galois extension $L/K$, there is an isomorphism $\inv_L : H^2(G_L, C^{\bullet}) \ra \Q/\Z$. We say that $(C^{\bullet}, \{\inv_L\})$ is a {\bf hyper-class formation} for $K$ if the following axioms are satisfied:

\begin{enumerate}
\item For every open subgroup $G_L$ of $G_K$, $H^1(G_L, C^{\bullet}) = 0$.
\item For every pair of open subgroups $G_{L} \subset G_{K'}$ of $G_K$, with $m= [L:K']$, the following diagram commutes:
$$
\begin{CD}
H^2(G_{K'}, C^{\bullet}) @>{\textup{Res}}>> H^2(G_{L}, C^{\bullet}) \\
@V{\inv_{K'}}VV @VV{\inv_{L}}V \\
\Q/\Z @>{\cdot m}>> \Q/\Z.
\end{CD}
$$
\end{enumerate}
\end{definition}


\begin{ex} 
Any class formation 
gives rise to a hyper-class formation, 
by considering
the 
formation module 
as a complex concentrated in degree zero.  There are three main examples of fields with class formations from the classical theory: local fields, global fields, and finite fields.  For a local field, the multiplicative group of the separable closure is a class formation; for a global field, the idele class group is a class formation; and, for a finite field, the trivial Galois module $\Z$ is a class formation.
\end{ex}

\begin{ex} \label{2loc} In \cite{Ko}, Koya gave the first 
new
example of hyper-class formation that is a genuine complex.  Let $K$ be a 2-local field, that is, a complete discrete valuation field with residue field a local field.  Then Lichtenbaum's weight two motivic complex $\Ze(2, K)$ is a hyper-class formation for $K$.  This complex was first constructed by Lichtenbaum in \cite{Li}; it is a two-term complex of $G_K$-modules, concentrated in degrees $-1$ and $0$.  Its homology is $K_2(K_{sep})$ in degree $0$, and $K_3^{ind}(K_{sep})$ in degree $-1$.  

The weight two motivic complex fits into a conjectural framework, see Lichtenbaum \cite{Li0}; in this setting, the weight one motivic complex $\Ze(1, K)$ is the multiplicative group of the separable closure, and $\Ze(0, K)$ is the trivial Galois module $\Z$, both considered as complexes concentrated in degree zero.  So, one should think of the hyper-class formation for a 2-local field as the natural extension of the class field theory for local and finite fields.
\end{ex}

\begin{ex}
More generally, if $K$ is an $n$-local field, that is, a complete discrete valuation field with residue field an $(n-1)$-local field, then the weight $n$ motivic complex 
is
a hyper-class formation for $K$, see Saito \cite{S4} and Spiess \cite{Sp}.
\end{ex}

If $C$ is a class formation, then $C^{G_L}$ is a class module for $\Gal(L/K)$, in the sense of Neukirch-Schmidt-Wingberg \cite{NSW}.  For complexes, we have the following analogous notion:

\begin{definition} Let $G$ be a finite group.  A {\bf class complex} for $G$ is a bounded 
below 
element in the derived category of 
$G$-modules $A^{\bullet}$ such that for all subgroups $H \subs G$, $H^1(H, A^{\bullet}) = 0$ and $H^2(H, A^{\bullet})$ is cyclic, of order $|H|$, generated by the restriction of a distinguished generator of $H^2(G, A^{\bullet})$.  
\end{definition}

\begin{remark}  If $C^{\bullet}$ is a hyper-class formation for $K$, and $L/K$ is any finite extension, then $\trg(G_L, C^{\bullet})$ is a class complex for $\Gal(L/K)$.  The proof, as in the corresponding classical statement, follows directly from the axioms using the inflation-restriction sequence of Galois cohomology. 

We will denote the canonical generator of $H^2(\Gal(L/K), \trg(G_L, C^{\bullet}))$ by $\al_{L/K}$, and refer to it as the fundamental class.  They satisfy the same compatibility conditions as their classical counterparts, see Serre \cite{Serre}.
\end{remark}

\begin{remark}  This operation for complexes is the analogue of taking Galois invariants: 
if $C^{\bullet}$ is concentrated in degree zero, then $\trg(G_L, C^{\bullet})$ is also a complex concentrated in degree zero, whose entry in degree zero is the $G_L$-invariants of the entry in degree 
zero
of $C^{\bullet}$.
\end{remark}

\begin{remark} \label{zar} 
This is closely related to Beilinson's motivic complexes over the Zariski site.  Let $K$ be a field, and $L/K$ a Galois extension.  Beilinson's weight two motivic complex $\Z_B(2, L)$ is a two-term complex of $\Gal(L/K)$-modules concentrated in degrees $-1$ and $0$.
Its homology is $K_2(L)$ in degree $0$, and $K_3^{ind}(L)$ in degree $-1$.
As noted by Lichtenbaum \cite{Li0}, the weight two Zariski motivic complex can be realized as the ``Galois invariants'' of the weight two $\e$tale motivic complex:
$$
\Z_B(2, L) = \trg(G_L, \Ze(2, K)).
$$
\end{remark}

\begin{ex}  Let $K$ be a complete discrete valuation field with residue field a global field, and let $L/K$ be a finite Galois extension.  
There is a form of class field theory for $K$: 
one can form an ``idele complex'' 
using a lifting construction of Kato \cite{K3} and the weight two Zariski motivic complex.  
This element of the derived category of $\Gal(L/K$-modules satisfies Galois descent, in the sense that $\trg(\Gal(L/K), I(L)^{\bullet})$ is quasi-isomorphic to $I(K)^{\bullet}$.  If $L/K$ is unramified, then the mapping cone of $\Z_B(2, L) \lra I(L)^{\bullet}$ is a formation complex for $L/K$, see \cite{me}. 
\end{ex}




\begin{remark} Classically, the Tate cohomology groups are at the heart of the 
cohomological interpretation of class field theory.  
We will see that a hyper-class formation gives rise to an Artin map using Tate cohomology of complexes.
\end{remark}

When $G$ is a finite group, we use the standard two-sided resolution $P_{\bullet}$ of $\Z$ as a $G$-module to define Tate cohomology groups, see Neukirch-Schmidt-Wingberg \cite{NSW}.   The following definition was 
made by Koya \cite{Ko}:

\begin{definition} Let $G$ be a finite group, $A^{\bullet}$ a bounded complex of $G$-modules, and $P_{\bullet}$ the standard two-sided resolution of $\Z$ as a $G$-module.  The Tate cohomology group $\HT^q(G, A^{\bullet})$ for $q \in \Z$ is defined to be the $q^{th}$ total homology of the double complex: 
$$
\Hom_G(P_{\bullet}, A^{\bullet}).
$$
\end{definition}

\begin{remark}  The assumption that $A^{\bullet}$ is a bounded complex, that is, that $A^{\bullet}$ has only finitely many non-zero entries, is necessary for the convergence of the spectral sequence associated to the double complex.  
If $A^{\bullet}$ is a complex that is zero in degree larger than $n$, then 
$$
\HT^q(G, A^{\bullet}) \simeq \Hc^q\big(\RG(G, A^{\bullet})\big) = H^q(G, A^{\bullet})
$$
for all $q > n$.  Furthermore, we have that 
$$
\HT^n(G, A^{\bullet}) = \Hc^n\big(\RG(G, A^{\bullet})\big)/N\Hc^n\big(A^{\bullet}\big), 
$$
see 
Theorem 1.2 of Koya \cite{Ko}.
\end{remark}

\begin{remark} 
Let $A$ be a $G$-module, and consider $A$ as a complex concentrated in degree zero.  Let $A[1]$ denote the shift of this complex down by one; then 
$$
\HT^{q}(G, A[1]) = \HT^{q+1}(G, A).
$$
Dimension-shifting 
implies that we
can produce a $G$-module $A'$ such that 
$$
\HT^q(G, A') \simeq \HT^q(G, A[1]).
$$  
We can construct 
such a modules by imbedding $A$ in an injective module $I$ and letting $A'$ be the cokernel.  In the language of complexes, 
we have that $A[1]$ is quasi-isomorphic to $I \ra A'$, considered as a complex concentrated in degrees $-1$ and $0$.  
We can compute the 
total
homology of the corresponding double complex by first taking vertical differentials.  But $I$ is injective, so $\HT^q(G, I) = 0$ for {\em all} $q \in \Z$.  Therefore, the spectral sequence degenerates after one step, giving that 
$$
\HT^q(G, A[1]) \simeq \HT^q(G, I \ra A') \simeq \HT^q(G, A').
$$
\end{remark}

For complexes in general, we have the following result:

\begin{prop}[Dimension-shifting for complexes] \label{cshift} Let $G$ be a finite 
group and $A^{\bullet}$ a bounded 
complex of $G$-modules.  Then there exists a $G$-module $A'$ such that for 
any subgroup $H$ of $G$, 
$$
\HT^q(H, A^{\bullet}) \simeq \HT^q(H, A')
$$
for {\em all} $q \in \Z$.
\end{prop}

\begin{proof} 
Choose an injective resolution $A^{\bullet} \ra I^{\bullet}$.  
Since $A^{\bullet}$ is bounded, 
we can choose $I^{\bullet}$ to 
be bounded below.  Suppose that $A^{\bullet}$ is zero above degree $n$, so that the truncated complex $\tau_{_{\leq n}} I^{\bullet}$ is 
quasi-isomorphic to $A^{\bullet}$.  For any subgroup $H$ of $G$, we 
have:
$$
\HT^q(H, A^{\bullet}) \simeq \HT^q(H, \tau_{_{\leq n}} I^{\bullet}).
$$

On the other hand, $\tau_{_{\leq n}} I^{\bullet}$ is a complex whose terms are all injective, except the $n^{th}$ term, call it $A$.  By definition, $\HT^q(H, \tau_{_{\leq n}} I^{\bullet})$ is the $q^{th}$ total homology of the double complex $\Hom_H(P_{\bullet}, \tau_{_{\leq n}} I^{\bullet})$.  
Taking 
vertical differentials, 
we see that 
the spectral sequence computing 
the
homology collapses after one step. 
Therefore, we have 
that
$$
\HT^q(H, A^{\bullet}) \simeq \HT^q(H, \tau_{_{\leq n}} I^{\bullet}) \simeq 
\HT^{q+n}(H, A).
$$ 
Then we can use standard dimension-shifting techniques for modules to 
find $A'$ as in the Proposition.
\end{proof}


In a certain sense, this explains why one can construct a Weil group associated to a hyper-class formation: at every finite level, we have a class {\em module} at our disposal, see Remark \ref{cmod}.  

One can also use dimension-shifting to reduce the following to the classical Tate-Nakayama Theorem:

\begin{thm}[Tate-Nakayama for complexes] \label{tn} Let $G$ be a finite group, and 
$$
A^{\bullet} \otimes B^{\bullet} \ra C^{\bullet}
$$
be a $G$-pairing of bounded complexes of $G$-modules, and let $\al$ be an element of $H^p(G, A^{\bullet})$.  Then for all subgroups $H \subs G$, cupping with $\Res_H(\al)$ induces a homomorphism
$$
c(\al)_{q,H} : \HT^q(H, B^{\bullet}) \ra \HT^{p+q}(H, C^{\bullet})
$$
for all $q \in \Z$.  Suppose that for some $q_0 \in \Z$ and every Sylow subgroup $H \subset G$, $c(\al)_{q_0 - 1, H}$ is surjective, $c(\al)_{q_0, H}$ is bijective, and $c(\al)_{q_0+1, H}$ is injective.  Then $c(\al)_{q, H}$ is an isomorphism for all $q \in \Z$ and all $H \subs G$.
\end{thm}

\begin{cor} \cite{Ko1}*{Thm. 6} \label{tnk}  Let $L/K$ be a finite Galois extension, and $A^{\bullet}$ a class complex for $\Gal(L/K)$, with $\al$ the distinguished generator of $H^2(\Gal(L/K), A^{\bullet})$.  Then
$$
\cdot \cup \Res_E(\al) : \HT^q(\Gal(L/E), \Z) \lra \HT^{q+2}(\Gal(L/E), A^{\bullet})
$$
is an isomorphism for all $q \in \Z$ and all intermediate fields $E$.
\end{cor}

As in the classical setting, Corollary \ref{tnk} allows one to construct 
an Artin map  associated to the hyper-class formation, as the inverse of 
a cup 
product on Tate cohomology.  
More precisely, if $\al_{L/K}$ is the fundamental class, we have an isomorphism:
$$
\cdot \cup \al_{L/K} : \HT^{-2}(\Gal(L/K), \Z) \lra \HT^0(\Gal(L/K), \trg(G_L, C^{\bullet})).
$$
The inverse to this cup product isomorphism gives rise to an ``Artin map'' in the limit:
$$
C(K) \lra G_K^{ab},
$$
where $C(K)$ denotes $H^0(G_K, C^{\bullet})$.

\begin{ex}
Let $K$ be a 2-local field, 
then the weight two motivic complex $\Ze(2, K)$ is a hyper-class formation, see Remark \ref{2loc}.  In this case, the Artin map provides a map between $G_K^{ab}$ and  
$$
H^0(G_K, \Ze(2, K)) = \Hc^0(\RG(G_K, \Ze(2, K))) = \Hc^0(\Z_B(2, K)),
$$
where $\Z_B(2, K)$ is the Zariski version of the weight two motivic complex.  The last equality holds because $\Z_B(2, K)$ can be realized as $\trg(G_K, \Ze(2, K))$, see Remark \ref{zar}.

But, we have that $\Hc^0(\Z_B(2, K)) = K_2(K)$, 
so that if $K$ is a 2-local field, then
there is an Artin map:
$$
K_2(K) \lra G_K^{ab},
$$
In fact, if $F$ is the residue field of $K$, then we have a commutative square:
$$ \begin{CD}
K_2(K) @>>> G_K^{ab} \\
@VVV @VVV \\
F^{\times} @>>> G_F^{ab},
\end{CD} $$
where the horizontal maps are Artin maps, and the first vertical map is the tame symbol, 
see Milnor \cite{Milnor}.  This diagram commutes because the fundamental classes for $K$ are compatible with the fundamental classes for the residue field $F$ under a boundary map that exists on the level of motivic complexes.
\end{ex}

\section{Weil groups} \label{weil}

In this Section, we first review the classical axiomatic construction of the Weil group associated to a field with a class formation.  Then, we show how one can also associate a Weil group to a field with a hyper-class formation.  

The primary object of study in this Section will be group extensions.  If $G$ is a finite group, and $A$ is a $G$-module, then an extension of $G$ by $A$ is an exact sequence:
$$
1 \ra A \lra U \lra G \ra 1.
$$
It is well-known that extensions of $G$ by $A$ up to equivalence are in one-to-one correspondence with classes in $H^2(G, A)$, see $\sect$ VII.3 of Serre \cite{Serre}.  

Let $K$ be a field, and $C$ a class formation for $K$.  Let $L/K$ be a finite Galois extension, and let $C(L)$ denote $C^{G_L}$.  Let $\al_{L/K}$ be the fundamental class, the canonical generator of $H^2(\Gal(L/K), C(L))$.  

\begin{definition} \label{wgpdef} A {\bf Weil group} $(\W(L/K), g, \{f_e\})$ for $L/K$ consists of:
\begin{enumerate}
\item A group $\W(L/K)$, referred to as the Weil group;
\item A homomorphism $g$ of $\W(L/K)$ onto $\Gal(L/K)$;
\item And, for each intermediate field $E$, an isomorphism 
$$
f_E : C(E) \approx \W(L/E)/\W(L/E)^c,
$$
where $\W(L/E)$ is defined to be $g^{-1}(\Gal(L/E))$;
\end{enumerate} 
such that,
\begin{enumerate}
\item For each intermediate extension $E'/E$, the following diagram commutes:
$$ \begin{CD} 
C(E) @>{f_E}>> \W(L/E)/\W(L/E)^c \\
@V{\incl}VV @VV{V}V \\
C(E') @>{f_{E'}}>> \W(L/E')/\W(L/E')^c,
\end{CD} $$
where $V$ is a transfer map;
\item For every $w \in \W(L/K)$, with $\sigma = g(w) \in \Gal(L/K)$, we have that $\W(L/E)^w$ is  $\W({L/E^{\sigma}})$ for every intermediate field $E$; and the following diagram commutes:
$$ \begin{CD}
C(E) @>{f_E}>> \W(L/E)/\W(L/E)^c \\
@V{\sigma}VV @VV{w}V \\
C(E^{\sigma}) @>{f_{E^{\sigma}}}>> \W(L/E^{\sigma})/\W(L/E^{\sigma})^c,
\end{CD} $$
where the second vertical map denotes conjugation by $w$;  
\item For any group $G$ with a subgroup $H$ of finite index, there is a natural group extension, called the factor subgroup extension, given by:
$$
1 \ra H/H^c \ra G/H^c \ra G/H \ra 1;
$$
for every intermediate Galois extension $E'/E$, we have that $\W(L/E')$ is a subgroup of finite index in $\W(L/E)$, with $\W(L/E)/\W(L/E')$ isomorphic to $\Gal(E'/E)$ via $g$; then the factor subgroup extension:
$$
1 \ra C(E') \ra \W(L/E)/\W(L/E')^c \ra \Gal(E'/E) \ra 1.
$$
corresponds to the fundamental class $\al_{E'/E}$; 
\item And, $\W(L/L)$ is abelian.  
\end{enumerate}
\end{definition}

\begin{remark} Applying the axioms, we see that 
$\W(L/K)$ is 
an extension of $\Gal(L/K)$ by $C(L)$, 
corresponding to the fundamental class $\al_{L/K}$.  \label{cfwgpdef} 
In fact, any extension 
$$ \begin{CD}
1 \lra C(L) @>>> \W(L/K) @>{g}>> \Gal(L/K) \lra 1 
\end{CD} $$
corresponding to $\al_{L/K}$ gives rise to a candidate for a Weil group.  
So, one should think of the Weil group as another manifestation of the fundamental class.

For any intermediate field $E$, we have the subextension 
$$
\W(L/E) = g^{-1}(\Gal(L/E))
$$
This group extension has corresponding cohomology class $\Res(\al_{L/K}) = \al_{L/E}$.  

Furthermore, for every intermediate field $E$, we have that the transfer map 
$$
V : \W(L/E)/\W(L/E)^c \lra C(L)
$$
induces an isomorphism of $\W(L/E)/\W(L/E)^c$ with $C(L)^{\Gal(L/E)} = C(E)$, see 
Chapter 12 of Artin-Tate \cite{AT}.  
A transfer map 
exists for any group $G$ with a subgroup $H$ of finite index; it is a map
$$
V : G/G^c \lra H/H^c
$$
that 
behaves like a norm map.  

\end{remark}

\begin{thm} \cite{AT}*{Ch. 14, Thm. 1} \label{cfwgp} Let $K$ be a field with a class formation $C$, 
and $L/K$ a finite Galois extension.  
There exists a Weil group $(\W(L/K), g, \{f_E\})$. 
\end{thm}

\begin{proof}  As in Remark \ref{cfwgpdef}, we choose an extension 
$$ \begin{CD}
1 \lra C(L) @>>> \W(L/K) @>{g}>> \Gal(L/K) \lra 1,
\end{CD} $$
whose corresponding cohomology class is the fundamental class $\al_{L/K}$.  Then we have isomorphisms
$$
f_E : C(E) \approx \W(L/E)/\W(L/E)^c,
$$
given as the inverses of the appropriate transfer maps.

Let $E'/E$ be an intermediate extension.  The first axiom follows 
from the 
commutative diagram:
$$ \begin{CD}
C(E) @<{V}<< \W(L/E)/\W(L/E)^c \\
@V{\incl}VV @VV{V}V \\
C(E') @<{V}<< \W(L/E')/\W(L/E')^c.
\end{CD} $$

Let $w \in \W(L/K)$, and put $\sigma = g(w)$.  Let $E$ be an intermediate field.  The second axiom follows 
from the 
commutative diagram:
$$ \begin{CD}
C(E)@<{V}<< \W(L/E)/\W(L/E)^c \\
@V{\sigma}VV @VV{w}V \\
C(E^{\sig}) @<{V}<< \W(L/E^{\sig})/\W(L/E^{\sig})^c.
\end{CD} $$

The third property follows from an analysis of factor subgroup extensions.  Let $\alpha \in H^2(\Gal(E'/E), C(E'))$ be the class corresponding to the extension:
$$
1 \ra C(E') \ra \W(L/E)/\W(L/E')^c \ra \Gal(E'/E) \ra 1.
$$
Then by Theorem 6 of Chapter 13 in Artin-Tate \cite{AT}, we have that 
$$
\infl(\al) = |E':E| \cdot \al_{L/E}.
$$
This is enough to imply that $\al$ is the fundamental class $\al_{E'/E}$.  

Finally, the fourth property to is trivial to verify; by construction, $\W(L/L)$ is $g^{-1}(1)$, which is 
$C(L)$.  
\end{proof}

Next, we present a construction of a Weil group associated to a field $K$ with a hyper-class formation $C^{\bullet}$.  Let $L/K$ be a finite Galois extension, and let $C(L)^{\bullet}$ denote $\trg(G_L, C^{\bullet})$.  Recall that $C(L)^{\bullet}$ is a class complex for $L/K$; the fundamental class $\al_{L/K}$ is the canonical generator of $H^2(\Gal(L/K), C(L)^{\bullet})$.  We have the standard spectral sequence:
$$
E_2^{p,q} = H^q(\Gal(L/K), \Hc^p(C(L)^{\bullet})) \Longrightarrow H^{p+q}(\Gal(L/K), C(L)^{\bullet}).
$$
Let $C(L)$ denote $\Hc^0(C(L)^{\bullet}) = H^0(G_L, C^{\bullet})$; 
there is an edge morphism:
$$
H^2(\Gal(L/K), C(L)^{\bullet}) \lra H^2(\Gal(L/K), C(L)),
$$
see Weibel \cite{W}.

\begin{definition} \label{hfun} The {\bf hyper-fundamental class} $\bg$ is defined to be the image of the fundamental class $\ag$ in $H^2(\Gal(L/K), C(L))$ under the edge morphism coming from the standard spectral sequence.  
\end{definition}

We proceed as before, and make the following construction:
\begin{definition} Let {\bf $\mathbf{\W(L/K)}$} be a group extension 
$$
\begin{CD}
1 \lra C(L) @>>> \W(L/K) @>g>> \Gal(L/K) \lra 1
\end{CD}
$$
whose corresponding cohomology class is the hyper-fundamental class $\bg$.  
\end{definition}

\begin{remark} \label{cmod}  Using dimension-shifting (Proposition 
\ref{cshift}), we have a $\Gal(L/K)$-module $A(L)$ such that
$$
\HT^q(\Gal(L/E), C(L)^{\bullet}) \simeq \HT^q(\Gal(L/E), A(L))
$$
for all $q \in \Z$ and all intermediate fields $E$.  
We recall the construction of $A(L)$: choose an injective resolution $C(L)^{\bullet} \ra I(L)^{\bullet}$ and let $A(L)$ be the entry in degree zero of $\tn I^{\bullet}$.  Note that $C(L)$ is naturally identified with the cokernel of the map $I(L)^{-1} \ra A(L)$,
and  
we can obtain the hyper-fundamental class $\bg$ as in the image of $\ag$ under the induced map 
$$
H^2(\Gal(L/K), A(L)) \ra H^2(\Gal(L/K), C(L)).
$$  
\end{remark}

\begin{lemma} Let $L/K$ be a finite Galois extension, $E$ an intermediate field, and $\W(L/K)$ a group extension corresponding to the hyper-fundamental class $\bg$.  Then the subgroup extension $g^{-1}(\Gal(L/E))$ has corresponding cohomology class $\bh$.
\end{lemma}

\begin{proof}  There is a map of group extensions:
$$
\begin{CD}
0 \lra \ @. C(L) @>>> g^{-1}(\Gal(L/E)) @>g>> \Gal(L/E) @. \ \lra 1 \\
@. || @. @VV{\textup{incl.}}V @VV{\textup{incl.}}V \\
0 \lra \ @. C(L) @>>> \W(L/K) @>g>> \Gal(L/K) @. \ \lra 1.
\end{CD}
$$
From this, we can conclude that the cohomology class of the top extension is the restriction of the class of the bottom extension.  But, 
$$
\ah = \res(\ag),
$$
which implies that the cohomology class corresponding to the top extension is the image of $\ah$.  
\end{proof}


The most basic property of a Weil group in the classical setting is that its abelianization is the module of the class formation.  In this setting, we have the following:

\begin{prop} Let $L/K$ be a finite Galois extension, and $E$ an intermediate field.  Then there is an isomorphism:
$$
f_E : C(E) \approx \W(L/E)/\W(L/E)^c,
$$
where $C(E) = H^0(G_E, C^{\bullet})$.
\end{prop}

\begin{proof}  
We can use our injective resolution $I(L)^{\bullet}$ to compute 
$$
C(E) = \Hc^0(\RG(\Gal(L/E), C(L)^{\bullet})). 
$$
We see that: 
$$
C(E) = \coker\Bigl[(I(L)^{-1})^{\Gal(L/E)} \ra A(L)^{\Gal(L/E)}\Bigr].
$$
Let $\U(L/E)$ be an extension of $\Gal(L/E)$ by $A(L)$ corresponding to $\al_{L/E}$.  Because $A(L)$ is a class module for $L/E$,   
we have that the transfer map 
$$
V : \U(L/E)/\U(L/E)^c \ra A(L)^{\Gal(L/E)}
$$
is an isomorphism, see Theorem 4 of Chapter 12 in Artin-Tate \cite{AT}.  Let $e_E$ be the inverse of this transfer map.  
We have a map of group extensions:
$$ \begin{CD}
1 \lra @. A(L) @>>> \U(L/E) @>{h}>> \Gal(L/E) @. \lra 1 \\
@. @VVV @VVV || \\
1 \lra @. C(L) @>>> \W(L/E) @>{g}>> \Gal(L/E) @. \lra 1,
\end{CD} $$
which induces a surjection:
$$
A(L)^{\Gal(L/E)} \approx \U(L/E)/\U(L/E)^c \ra \W(L/E)/\W(L/E)^c.
$$
To prove the Proposition, we will show that the kernel of this map can be identified with $(I(L)^{-1})^{\Gal(L/E)}$.



By a simple diagram chase, we can identify the kernel with:
\begin{eqnarray*}
N\Bigl(\image\bigl[I(L)^{-1} \ra A(L)\bigr]\Bigr) &=& \image\Bigl[N\bigl(I(L)^{-1}\bigr) \ra A(L)^{\Gal(L/E)}\Bigr] \\
&=& \image\Bigl[(I(L)^{-1})^{\Gal(L/E)} \ra A(L)^{\Gal(L/E)}\Bigr], \\
\end{eqnarray*}
where the last equality follows because $I(L)^{-1}$ is injective. 
So, $e_E$ gives rise to 
the desired isomorphism
$$
f_E  :  C(E) \approx \W(L/E)/\W(L/E)^c.
$$
\end{proof}


\begin{thm} \label{hweil}  The collection $(\W(L/K), g, \{f_E\})$ associated to a field $K$ with a hyper-class formation $C^{\bullet}$ satisfies the axioms for a Weil group.
\end{thm}

\begin{proof}  Verification of the axioms follows along the lines of Theorem \ref{cfwgp}.  The first two axioms follow from naturality of the transfer map, the third follows from some compatibility between hyper-fundamental classes, and the last is trivial.


Let $E/E'$ be an intermediate extension.  The first axiom follows from the commutative diagram:
$$ \begin{CD}
A(L)^{\Gal(L/E)} @<{V}<< \U(L/E)/\U(L/E)^c @. \ \lra \ @. \W(L/E)/\W(L/E)^c \\
@V{\textup{incl}}VV @VV{V}V @. @VV{V}V \\
A(L)^{\Gal(L/E')} @<{V}<< \U(L/E')/\U(L/E')^c @. \ \lra \ @. \W(L/E')/\W(L/E)^c.
\end{CD} $$


Let $w \in \W(L/K)$, and put $\sig = g(w)$.  Let $u \in \U(L/K)$ be any lift of $w$, and let $E$ be an intermediate field.  The second axioms follows from the commutative squares:
$$ \begin{CD}
A(L)^{\Gal(L/E)} @<{V}<< \U(L/E)/\U(L/E)^c \\
@V{\sigma}VV @VV{u}V \\
A(L)^{\Gal(L/E^{\sigma})} @<{V}<< \U(L/E^{\sigma})/\U(L/E^{\sigma}),
\end{CD} $$
and
$$ \begin{CD}
\U(L/E)/\U(L/E)^c @>>> \W(L/E)/\W(L/E)^c \\
@VV{u}V  @VV{w}V \\
\U(L/E^{\sigma})/\U(L/E^{\sigma}) @>>> \W(L/E^{\sigma})/\W(L/E^{\sigma})^c.
\end{CD} $$

The third axiom is perhaps the most interesting, because implies that a Weil group for $L/K$ gives rise to a Weil group for any intermediate extension $E'/E$.  
It can be thought of as a ``compatibility condition'' for the hyper-fundamental classes.  However, 
in general there is no natural map between $H^2(\Gal(E'/E), C(E))$ and $H^2(\Gal(L/E), C(L))$.  

Let $\beta$ denote the cohomology class of the factor subgroup extension:
$$
1 \ra C(E') \ra \W(L/E)/\W(L/E')^c \ra \Gal(E'/E) \ra 1
$$
We have a map of group extensions:
$$ \begin{CD}
1  \ra \ @.  A(L)^{\Gal(L/E')} @. \ \lra \ @. \U(L/E)/\U(L/E')^c @. \ \lra \ @. \Gal(E'/E) @. \ \ra 1 \\
@. @VVV @. @VVV @. || \\
1 \ra \ @. C(E') @. \lra \ @. \W(L/E)/\W(L/E')^c @. \ \lra \ @. \Gal(E'/E) @. \ \ra 1.
\end{CD} $$
Let $\alpha$ 
denote the cohomology class of the top extension.  By generalities on factor subgroup extensions, we have that 
$$
\infl(\alpha) = [L:E'] \cdot \alpha_{L/E}.
$$

On the other hand, we have a commutative diagram:
$$\begin{CD} H^2(\Gal(E'/E), C(E')^{\bullet}) @>{\textup{Infl}}>> H^2(\Gal(L/E), C(L)^{\bullet}) \\
@AAA @AAA \\
H^2(\Gal(E'/E), A(L)^{\Gal(L/E')}) @>{\textup{infl}}>> H^2(\Gal(L/E), A(L)), 
\end{CD} $$
where the vertical maps are induced by the natural map
$$
A(L) \lra C(L)^{\bullet}.
$$
Now, the axioms for a hyper-class formation imply that the top inflation map is multiplication by $[L:E']$.  That is, the class $\alpha$ must be mapped to the fundamental class $\alpha_{E'/E} \in H^2(\Gal(E'/E), C(E')^{\bullet})$.  

But the natural map 
$$
H^2(\Gal(E'/E), A(L)^{\Gal(L/E')}) \ra H^2(\Gal(E'/E), C(E'))
$$
coming from the surjection $A(L)^{\Gal(L/E')} \ra C(E')$ factors as a composition:
\begin{eqnarray*}
H^2(\Gal(E'/E), A(L)^{\Gal(L/E')}) &\lra& H^2(\Gal(E'/E), C(E')^{\bullet}) \\ 
&\lra& H^2(\Gal(E'/E), C(E'))
\end{eqnarray*}
where is second map is an edge morphism.  This means precisely that $\beta$ is the image of the fundamental class $\alpha_{E'/E}$; that is, that $\beta$ is the hyper-fundamental class $\beta_{E'/E}$.

Finally, 
the fourth axiom is trivial 
by construction $\W(L/L) = g^{-1}(1) = C(L)$ 
is abelian.  
\end{proof}

\begin{remark} Because the Weil group of a hyper-class formation satisfies the classical axioms, we can conclude that it is unique up to Weil isomorphism, see Chapter 14 of Artin-Tate \cite{AT}.  
In particular, 
our construction does not depend in an essential way on the choice of injective resolution.
\end{remark}

\begin{remark} The most remarkable property of the Weil group is that it captures the Artin map in a natural way, see Theorem 5 of Chapter 14 in Artin-Tate \cite{AT}.  
If $L/K$ is a finite Galois extension, then the Artin map for $L/K$ is given by:
$$ \begin{CD}
C(K)/NC(L) @>{f_K}>> \W(L/K)/(\W(L/L) \cdot \W(L/K)^c) \approx \Gal(L/K)^{ab}. 
\end{CD} $$
\end{remark}

\section{Examples}

The purpose of this Section is to describe the Weil groups associated to the two examples of hyper-class formations, and their connection to the classical examples of Weil groups.

\begin{ex} Let $l/k$ be a finite Galois extension of finite fields.  
The Weil group $\W(l/k)$ is completely determined by the degree $n=[l:k]$.  In this case, $\Gal(l/k) \simeq \frac{1}{n}\Z/\Z$, 
and 
$\W(l/k)$ 
is an 
extension of $\frac{1}{n}\Z/\Z$ by $\Z$.  Because the Galois group is abelian, and it acts trivially on $\Z$, it follows directly from the axioms that the Weil group for $l/k$ is the extension:
$$
0 \lra \Z \lra \frac{1}{n}\Z \lra \frac{1}{n}\Z/\Z \lra 0.
$$
\end{ex}

\begin{ex} Let $E/F$ be a finite Galois extension of local fields, with residue field extension $l/k$.  Then $\W(E/F)$ will be a group extension of $\Gal(E/F)$ by $E^{\times}$.  There is a map of group extensions:
$$ \begin{CD}
1 @>>> E^{\times} @>>> \W(E/F) @>>> \Gal(E/F) @>>> 1 \\
 @. @V{v_E}VV @VVV @VVV \\
0 @>>> \Z @>>> \W(l/k) @>>> \Gal(l/k) @>>> 1.
\end{CD} $$

For $E/F$ unramified, this follows directly from the proof of local class field theory given in $\sect$ XIII.3 of Serre \cite{Serre}: the isomorphism 
$$ \begin{CD}
H^2(\Gal(E/F), E^{\times}) @>{{v_E}_*}>> H^2(\Gal(l/k), \Z)
\end{CD} $$
is used to define the fundamental class $\al_{E/F}$.  Using properties of group extensions, the fact that $\al_{E/F}$ maps to $\al_{l/k}$ under ${v_E}_*$ is enough to imply the existence a map from $\W(E/F)$ to $\W(l/k)$ making the diagram commute.

If $E/F$ is any finite Galois extension, we let $E'$ be the maximal unramified extension of $F$ contained in $E$.  Then there is a commutative diagram:
$$ \begin{CD}
1 @>>> E'^{\times} @>>> \W(E'/F) @>>> \Gal(E'/F) @>>> 1 \\
@. @V{v_{E'}}VV @VVV || \\
0 @>>> \Z @>>> \W(l/k) @>>> \Gal(l/k) @>>> 1.
\end{CD} $$
But, we can use the third axiom to relate $\W(E'/F)$ with $\W(E/F)$.  More precisely, we can obtain $\W(E'/F)$ from $\W(E/F)$ as a factor subgroup extension.  This amounts to a commutative diagram:
$$ \begin{CD}
1 @>>> E^{\times} @>>> \W(E/F) @>>> \Gal(E/F) @>>> 0 \\
@. @V{N_{E/E'}}VV @VVV @VVV \\
1 @>>> E'^{\times} @>>> \W(E'/F) @>>> \Gal(E'/F) @>>> 0.
\end{CD} $$
Then the original claim follows by putting together the diagrams, along with the well-known fact that $v_{E} = v_{E'} \circ N_{E/E'}$ (see $\sect$ II.2 of Serre \cite{Serre}).
\end{ex}

\begin{ex} Let $K$ be a 2-local field, that is, a complete discrete valuation field with residue field a local field, and let $L/K$ be a finite Galois extension.  Then the weight two Zariski motivic complex $\Z_B(2, L)$ (see Remark \ref{zar}) is a formation complex for $L/K$.  So, the associated Weil group $\W(L/K)$ will be a group extension of $\Gal(L/K)$ by $\Hc^0(\Z_B(2, L)) = K_2(L)$.  Let $E/F$ be the residue field extension, and assume that this extension is separable.  
As before, there is a map of group extensions:
$$ \begin{CD}
1 @>>> K_2(L) @>>> \W(L/K) @>>> \Gal(L/K) @>>> 1 \\
@. @V{\delta_L}VV @VVV @VVV \\
1 @>>> E^{\times} @>>> \W(E/F) @>>> \Gal(E/F) @>>> 1,
\end{CD} $$
where $\delta_L$ denotes the tame symbol.

This can be proven using the method of the previous example.  The unramified case follows from 
the 
proof of class field theory for 2-local fields.  
The general case follows from the axioms, and the fact that if $L/L'$ is a totally ramified extension, then $\delta_L = \delta_{L'} \circ N_{L/L'}$, see Lemma 16 of $\sect$ 1.7 in Kato \cite{K2}.
\end{ex}

\begin{ex} Let $K$ be a complete discrete valuation field with residue field a global field, and let $L/K$ be a finite {\em unramified} extension with residue field extension $E/F$.  Then the ``idele class complex'' $C(L)^{\bullet}$ is a formation complex for $L/K$, and the associated Weil group $\W(L/K)$ will be a group extension of $\Gal(L/K)$ by Kato's $K_2$-idele class group $C_2(L)$, see \cite{K3}.  

It follows directly from the computations in \cite{me} that there is a map of group extensions:
$$ \begin{CD}
1 @>>> C_2(L) @>>> \W(L/K) @>>> \Gal(L/K) @>>> 1 \\
@. @V{\prod_{\om} \delta_{\om}}VV @VVV || \\
1 @>>> C(E) @>>> \W(E/F) @>>> \Gal(E/F) @>>> 1,
\end{CD} $$
where $C(E)$ is the usual idele class group associated to a global field $E$, and $\delta_{\om}$ denotes the local tame symbol.
\end{ex}

\end{document}